\documentclass[12pt]{article}

\usepackage{amsthm,amssymb,amsfonts,latexsym}
\newtheorem{theorem}{Theorem}

\begin{document}

\overfullrule=0pt
\baselineskip=24pt
\font\tfont= cmbx10 scaled \magstep3
\font\sfont= cmbx10 scaled \magstep2
\font\afont= cmcsc10 scaled \magstep2
\title{\tfont On relations of invariants for 
vector-valued forms }
\bigskip
\author{ Thomas Garrity \\
Zachary Grossman\footnote{Current address: Department of Economics, 
University of California, Berkeley}
\\ Department of Mathematics\\ Williams College\\Williamstown, MA  
01267\\ 
email:tgarrity@williams.edu}
\date{}
\maketitle
\begin{abstract}
An algorithm is given for computing explicit formulas for the 
generators of relations among the invariant rational functions  for vector-valued bilinear forms.
  These formulas have applications in the geometry of 
Riemannian submanifolds and in CR geometry.

\end{abstract}

\section{Preliminaries}
\subsection{Introduction}

Vector-valued forms play a key role in the study of higher 
codimensional geometries.  For example, they
 occur naturally in the study of 
Riemannian submanifolds (as the second fundamental form) and in CR 
geometry (as the Levi form).  In each of these there are natural 
group actions acting on the vector-valued forms,
 taking care of different choices of local coordinates 
and the such.  The algebraic invariants of these forms under these group 
actions provide invariants for the given geometries.  In
Riemannian geometry, for example, the scalar curvature can be expressed as an 
algebraic invariant of the second fundamental form.
But before these  invariants can be 
used, their algebraic structure must be known.  In \cite{GM}, an explicit list 
of the generators is given.  In that paper, though, there is no hint as 
to the relations among these generators.  In this 
paper, a method is given for producing  a list of the generators for the 
relations of the invariants.

In \cite{GM}, the problem of finding the rational invariants of bilinear maps 
from a complex vector space $V$ of dimension $n$ to a complex
vector space $W$ of dimension $k$, on which the group 
$GL(n,\mbox{{\bf C}})\times GL(k, \mbox{{\bf C}})$ acts, is reduced 
 to the problem of 
finding invariant one-dimensional subspaces of the vector spaces
$(V\otimes V\otimes W^{*})^{\otimes r}$, for  each positive integer   $r$.
  From this, it is shown that the 
invariants can be interpreted as being  generated by 
$$(\mbox{Invariants for $GL(n, \mbox{{\bf C}})$ of $V\otimes V$})\otimes
(\mbox{Invariants for $GL(k,\mbox{{\bf C}})$ of $W^{*}$}),$$
each component of which had been computed classically. 

 In this paper 
we extend this type of result, showing how to compute the relations
of bilinear forms from knowledge of
 the relations for $V \otimes V$ under the 
action of 
$GL(n, \mbox{{\bf C}})$ and the relations for $W^{*}$ under the 
action of 
$GL(k, \mbox{{\bf C}})$.  In particular, in a way that will be made 
more precise later, we show that the relations can be interpreted as 
being 
generated by:
$$((\mbox{Relations for}\;GL(n, \mbox{\bf {C}})\; \mbox{of}\;
V\otimes V)\otimes
(\mbox{Generators for}\; GL(k, \mbox{\bf {C}}) \;\mbox{of}\;W^{*}))$$
$$\oplus ((\mbox{Generators for}\; GL(n, \mbox{\bf {C}})\; \mbox{of}\;
V\otimes V)\otimes
(\mbox{Relations for}\; GL(k, \mbox{\bf {C}}) \;\mbox{of}\;W^{*})),$$
each component of which is  known classically.  
While this paper concentrates on the case when $GL(n, \mbox{\bf {C}})$ acts on the
vector space $V$ and $GL(k, \mbox{\bf {C}})$ acts on $W$, the techniques that we use
are  applicable for when $G$ and $H$ are any completely reducible Lie groups
acting on the vector spaces $V$ and $W$, respectively.  When the bilinear form
is the second fundamental form of Riemannian geometry, then $G$ is the orthogonal
group $O(n)$ and $H$ is the orthogonal group $O(k)$.  In CR geometry, when the bilinear
form is the Levi form, then $G= GL(n, \mbox{\bf {C}})$, but now $H=GL(k, \mbox{\bf {R}})$.

In sections 1.2 through 1.4, we set up our basic notation.  Section 2 recalls how
to compute the invariants of bilinear forms.  Section 3 recalls the classically known 
relations for invariants of the general linear group.  While well-known, we spend time
on writing these relations both in the bracket notation for vectors and in the tensor
language that we are interested in.  Section 4 gives the relations 
among the invariant one dimensional subspaces of $(V\otimes V\otimes 
W^{*})^{\otimes r}$, for each positive integer $r$.  This is the key 
step in this paper. The key proof will be seen to be not hard, 
reflecting the fact that the difficulty in this paper is not the 
proofs but the finding of the correct statements and correct 
formulations of the theorems.  Section 5  gives a concrete example of a relation 
from section 4.  Section 6 finally deals with the finding the 
relations of the invariants for vector-valued bilinear forms.  
Section 7 gives concrete, if not painful, examples and discusses some 
geometric insights behind some of the computed invariants and relations.
Section 8 closes with some further questions.

It appears that the closest earlier work to this paper is in the study 
of the invariants of $n\times n$ matrices (see \cite{Formanek} and 
\cite{Procesi}), but the group actions are 
different in this case and links are not apparent. For general 
background in invariant theory, see \cite{DC}, \cite{Popov}, \cite {Olver} or 
\cite{Weyl}.

The authors would like to thank the referees for major help in improving the exposition
of this paper.  T. Garrity would also like to thank the mathematics department of
the University of Michigan, where part of this paper was written during a delightful sabbatical.

\subsection{Vector-valued forms}
For the rest of this paper, let $V$ be a complex $n$-dimensional vector space
  and $W$ be a complex $k$-dimensional vector space.  We are concerned 
  with the vector space $Bil(V,W)$, the space of bilinear maps from 
  $V\times V$ to $W$.  Each such bilinear map is an element of 
  $V^{*}\otimes V^{*}\otimes W$, where $V^{*}$ is the dual space of 
  $V$. The group $Aut(V) \times Aut(W)$ acts on 
  $Bil(V,W)$ by 
  $$gb(x,y)=pb(a^{-1}x,a^{-1}y)$$ for all
$g=(a,p)\in Aut(V) \times Aut(W)$, $b\in Bil(V,W)$, and $x,y \in V$.
Stated differently,  we  define
 $gb$ so that the following diagram commutes: $$a \> \times a \>
\begin{array}{ccc}
 V\times V & \stackrel{b}{\rightarrow} & W \\
 \downarrow &  & \downarrow \\
 V\times V & \stackrel{\rightarrow}{gb} & W
\end{array} p $$

As mentioned in the introduction, the results of this paper (and the results in  \cite{GM})
 also work for completely reducible 
subgroups of $Aut(V)$ and 
$Aut(W)$, though for 
simplicity, we  restrict our attention to the full groups $Aut(V)$ 
and $Aut(W)$, which of course are isomorphic to $GL(n,\mbox{\bf {C}})$ and 
$GL(k,\mbox{\bf {C}})$.


\subsection{Invariants}

Let $G$ be a group that acts linearly on a complex vector space $V$.  
A function
$$f:V \rightarrow \mathbb{C}$$ 
is a (relative) invariant if for all $g \in G$ and all $v\in V$, we 
have
$$f(gv) = \chi (g)f(v),$$
where $\chi:G \rightarrow \mathbb{C} - \{ 0 \}$ is a homomorphism 
(i.e. $\chi$ is an abelian character for the group $G$).  We call 
$\chi$ the {\it weight} of the invariant.  Note that the sum of two 
invariants of the same weight is another invariant.  Thus the 
invariants of the same weight will form in natural way a vector space. 

 As seen in 
\cite{DC} on pp. 5-9,
every rational invariant is the quotient of polynomial invariants, 
every polynomial invariant is the sum of homogeneous polynomial 
invariants, every degree $r$ homogeneous polynomial corresponds to an 
invariant $r$-linear function on the Cartesian product $V^{\times r}$ 
and every invariant $r$-linear function on  $V^{\times r}$ 
corresponds to an invariant linear function on  the $r$-fold 
tensor product $V^{\otimes r}$.  Thus to study rational invariants 
on $V$ we 
can concentrate on understanding the invariant one-dimensional 
subspaces on $V^{*\otimes r}$.

Let $\mathbb{C}[V]$ be the algebra of polynomial functions on $V$ and let
 $\mathbb{C}[V]^{G}$
denote the  algebra of the polynomials invariant under the action 
of $G$.  In general, the goal of invariant theory is
 to find a list of generators of  
algebra $\mathbb{C}[V]^{G}$ (a ``First Fundamental Theorem"), a list 
of generators of the relations of these generators (a ``Second 
Fundamental Theorem") and then relations of relations, etc.  A full 
such description is the   syzygy of $\mathbb{C}[V]^{G}$, etc.  

 By 
the above, we need to find the homogeneous  polynoials in 
$\mathbb{C}[V]^{G}$.  Since the homogeneous polynomials of degree r in $\mathbb{C}[V]$
are isomorphic to symmetric tensors in $V^{*\odot r}$,  we need to 
find the invariant one-dimensional 
subspaces of $V^{*\odot r}$.  Now, if $G$ acts on $V$, it will act on $V^{*\otimes r}$.
Suppose we have all invariant one-dimensional 
subspaces on $V^{*\otimes r}$.  Then we can easily recover all invariant one-dimensional 
subspaces on $V^{*\odot r}$ by the symmetrizing map from $V^{*\otimes 
r}$ to $V^{*\odot r}$.   This is the procedure we will use.

We are interested in rational invariants for the vector space $Bil(V,W)$ 
under the group action of $Aut(V) \times Aut(W)$. Thus we are interested 
in rational invariants for the vector space $Bil(V,W)$.  Hence the 
vector space we are interested in is $V^{*}\otimes V^{*}\otimes W$ 
under the group action of $Aut(V) \times Aut(W)$, which we will see 
means that we are interested initially in the  
invariant one-dimensional subspaces of $(V\otimes V \otimes 
W^{*})^{\otimes r}$ for each $r$  and 
finally  in the invariant one-dimensional subspaces of $(V\otimes V \otimes 
W^{*})^{\odot r}$.

\subsection{Indicial notations}

The following permutation notation will be used heavily throughout
this paper.  For any positive
integer $m$, define the permutation symbol
$\varepsilon^{i_{1}i_{2}\ldots i_{m}}$ to be equal to 1 if
$i_{1}i_{2}\ldots i_{m}$ is an even permutation of $1,2,\ldots ,m$,
to be equal to $(-1)$ if it is an odd permutation, and to be equal
to 0 otherwise. To indicate the product of $d$ (where $d$ is a
positive integer) such symbols for any permutation $\sigma \in
S_{dm}$, we use the shorthand notation $$\varepsilon^{I}
(m,dm,\sigma) =\varepsilon^{i_{\sigma(1)}\dots i_{\sigma(m)}}
\varepsilon^{i_{\sigma(m+1)}\dots i_{\sigma(2m)}}\cdots
\varepsilon^{i_{\sigma(dm-m+1)}\dots i_{\sigma(dm)}}.$$ 
The symbols $\varepsilon_{i_{1}i_{2}\dots i_{m}}$ and 
$\varepsilon_{I}
(m,dm,\sigma)$ are defined in a similar manner.
 The Einstein summation notation will be used.  Thus whenever a 
superscript and a subscript appear in the same term, this means sum 
over that term. 

As an example of the notation,
let $V$ be a two dimensional vector space with the basis $\{
e_{1},e_{2}\}$. Then $\varepsilon^{I} (2,2,\mbox{identity}
)e_{i_{1}} \otimes e_{i_{2}}$ denotes the following summation of
two-tensors from $V\otimes V$:
\begin{eqnarray*}
\varepsilon^{I} (2,2,\mbox{identity} )e_{i_{1}} \otimes e_{i_{2}}
& = & \; \varepsilon^{11}e_{1}\otimes e_{1} +\varepsilon^{12}
e_{1} \otimes e_{2} \\ &  & +\varepsilon^{21}e_{1}\otimes e_{2}
+\varepsilon^{22} e_{2} \otimes e_{2} \\
 & = & e_{1}\otimes e_{2} - e_{2}\otimes e_{1} 
 \\ & = &  e_{1}\wedge  e_{2}
\end{eqnarray*}

A slightly more complicated example is $\varepsilon^{I} (2,4,\mbox{identity}
)e_{i_{1}} \otimes e_{i_{2}} \otimes e_{i_{3}} \otimes e_{i_{4}}$.  
In  $\varepsilon^{i_{1}i_{2}, i_{3}, i_{4}}$, each $i_{m}$ can be 
either a $1$ or a $2$.  Thus there are $2^{4}$ terms being summed.  
But whenever at least three of the $i_{m}$ are a $1$ or $2$, the 
corresponding term is zero.  Hence there are really only six terms 
making up $\varepsilon^{I} (2,4,\mbox{identity}
)e_{i_{1}} \otimes e_{i_{2}} \otimes e_{i_{3}} \otimes e_{i_{4}}$.    We have
\begin{eqnarray*}
\varepsilon^{I} (2,4,\mbox{identity} )e_{i_{1}} \otimes e_{i_{2}}\otimes e_{i_{3}} \otimes e_{i_{4}}
& = & \; \varepsilon^{11}\varepsilon^{22}e_{1}\otimes e_{1}\otimes 
e_{2}\otimes e_{2}  \\
& & +\:\varepsilon^{12}\varepsilon^{12}
e_{1} \otimes e_{2} \otimes e_{1} \otimes e_{2} \\ 
&  & +\: \varepsilon^{12}  \varepsilon^{21}e_{1}\otimes e_{2}\otimes 
e_{2}\otimes e_{1} \\
& & +\: \varepsilon^{21}  \varepsilon^{12}e_{2}\otimes e_{1}\otimes 
e_{1}\otimes e_{2} \\
& & + \:\varepsilon^{21}  \varepsilon^{21}e_{2}\otimes e_{1}\otimes 
e_{2}\otimes e_{1} \\
& & +\:\varepsilon^{22}  \varepsilon^{11}e_{2}\otimes e_{2}\otimes 
e_{1}\otimes e_{1} \\
& =&e_{1} \otimes e_{2} \otimes e_{1} \otimes e_{2} \\ 
&  & -\:e_{1}\otimes e_{2}\otimes 
e_{2}\otimes e_{1} \\
& & -\:e_{2}\otimes e_{1}\otimes 
e_{1}\otimes e_{2} \\
& & +\: e_{2}\otimes e_{1}\otimes 
e_{2}\otimes e_{1} \\
&=&  (e_{1}\wedge  e_{2})\otimes (e_{1}\wedge e_{2}).
\end{eqnarray*}
We will see in  section three that this will be the invariant on four 
vectors ${\bf v}_1,{\bf v}_2,{\bf v}_3,{\bf v}_4$ corresponding to the product of 
determinants:
$$[{\bf v}_1, {\bf v}_2][{\bf v}_3,{\bf v}_4].$$

One more example that we will use later.  Consider $\varepsilon^{I} (2,4,(23)
)e_{i_{1}} \otimes e_{i_{2}} \otimes e_{i_{3}} \otimes e_{i_{4}}$. 
All we need to do is to flip, in the above formulas, $i_{2}$ with $i_{3}$ in the $\varepsilon^{i_{1}i_{2}, i_{3}, i_{4}}$.
Thus

\begin{eqnarray*}
\varepsilon^{I} (2,4,(23 )e_{i_{1}} \otimes e_{i_{2}}\otimes e_{i_{3}} \otimes e_{i_{4}}
& = & \; \varepsilon^{12}\varepsilon^{12}e_{1}\otimes e_{1}\otimes 
e_{2}\otimes e_{2}  \\
& & +\:\varepsilon^{11}\varepsilon^{22}
e_{1} \otimes e_{2} \otimes e_{1} \otimes e_{2} \\ 
&  & +\: \varepsilon^{12}  \varepsilon^{21}e_{1}\otimes e_{2}\otimes 
e_{2}\otimes e_{1} \\
& & +\:\varepsilon^{21}  \varepsilon^{12}e_{2}\otimes e_{1}\otimes 
e_{1}\otimes e_{2} \\
& & + \:\varepsilon^{22}  \varepsilon^{11}e_{2}\otimes e_{1}\otimes 
e_{2}\otimes e_{1} \\
& & +\:\varepsilon^{21}  \varepsilon^{21}e_{2}\otimes e_{2}\otimes 
e_{1}\otimes e_{1} \\
& =& e_{1} \otimes e_{1} \otimes e_{2} \otimes e_{2} \\ 
&  & -\: e_{1}\otimes e_{2}\otimes 
e_{2}\otimes e_{1} \\
& & -\:e_{2}\otimes e_{1}\otimes 
e_{1}\otimes e_{2} \\
& & + \: e_{2}\otimes e_{2}\otimes 
e_{1}\otimes e_{1} .
\end{eqnarray*}

In    section three we will see that this will be the invariant on  
the four 
vectors ${\bf v}_1,{\bf v}_2,{\bf v}_3,{\bf v}_4$ corresponding to the product of 
determinants:
$$[{\bf v}_1, {\bf v}_3][{\bf v}_2,{\bf v}_4].$$

\section{Generators for invariants of $Bil(V,W)$ }

This section is a quick review of the notation and the results in \cite{GM}, which we
need for the rest of this paper.

Let $e_{1},\ldots,e_{n}$ and $f_{1},\ldots, f_{k}$ be bases for $V$ 
and $W$ and $e^{1},\ldots,e^{n}$ and $f^{1},\ldots, f^{k}$
be dual bases for $V^{*}$ and $W^{*}.$   The  goal in \cite{GM} is to find the 
invariant one-dimensional subspaces of $(V\otimes V\otimes 
W^{*})^{\otimes r}$, for each possible $r$. We will throughout 
regularly identify $(V\otimes V\otimes 
W^{*})^{\otimes r}$ with $V^{\otimes 2r}\otimes (W^{*})^{\otimes r}$.

Let $r$ be a positive integer such that $n$ divides $2r$ and $k$ 
divides $r$.  For any element $\sigma$ in the permutation group $S_{2r}$ 
and any element $\eta$ in $S_{r}$, define
$$v_{\sigma} = \varepsilon^{I}(n, 2r, \sigma) e_{i_{1}}\otimes \ldots 
\otimes e_{i_{2r}}$$
and 
$$w^{\eta}= \varepsilon_{J}(k, r, \eta) f^{i_{1}}\otimes \ldots 
\otimes f^{i_{r}}.$$

\begin{theorem}
The vector space $ V^{\otimes 2r}\otimes (W^{*})^{\otimes r} $ 
has an invariant one-dimensional subspace if and 
only if $n$ divides $2r$ and $k$ 
divides $r$.  Every invariant one-dimensional subspace is a linear 
combination of various $v_{\sigma}\otimes w^{\eta}$, where $\sigma$ 
and $\eta$ range   through $S_{2r}$ and $S_{r}$, respectively.  
\end{theorem}

For each $r$,  denote the subspace generated by all of the various 
$v_{\sigma}\otimes w^{\eta}$ in $ V^{\otimes 2r}\otimes (W^{*})^{\otimes r} $ 
by 

$$( V^{\otimes 2r}\otimes (W^{*})^{\otimes r}  )_{\mbox{inv}}.$$
As shown in \cite{GM}, for each r the corresponding weights are the 
same.  Hence the sum of any two $v_{\sigma}\otimes w^{\eta}$ spans 
another one-dimensional
invariant subspace of $ V^{\otimes 2r}\otimes (W^{*})^{\otimes r} $ .  
Thus, for each $r$,  $( V^{\otimes 2r}\otimes (W^{*})^{\otimes r}  )_{\mbox{inv}}$
is the invariant subspace of $ V^{\otimes 2r}\otimes (W^{*})^{\otimes r} $ 
under our group action.

Hence for each r, the theorem is giving us a spanning set 
for  $( V^{\otimes 2r}\otimes (W^{*})^{\otimes r}  )_{\mbox{inv}}$.  Part of the goal 
of this paper is to produce an algorithm to find the relations among 
the elements for these spanning sets.

Let us put this into the language of bilinear forms, which will aid us 
later we look at specific examples.  By making our 
choice of bases, we can write each bilinear map from $V\times V$ to 
$W$ as a  $k$-tuple of $n\times n$ matrices $(B^{1},\ldots,B^{k})$, 
where each $B^{\alpha}= (B^{\alpha}_{ij})$.  More precisely, if $b\in 
Bil(V,W)$, then $ B^{\alpha}_{ij} = f^{\alpha}b(e_{i},e_{j})$.
We can restate the above theorem in terms of the $B^{\alpha}_{ij}$.

\begin{theorem}
There exists  a nonzero homogeneous invariant of degree $r$ on 
$Bil(V,W)$ only if $n$ divides $2r$ and $k$ 
divides $r$.  Further, every such  homogeneous invariant is 
a linear combination of various $f^{\sigma}_{\eta}$, where
$$f^{\sigma}_{\eta}= \varepsilon^{I}(n, 2r, \sigma) \varepsilon_{J}(c, r, \eta)
B^{j_{1}}_{i_{1}i_{2}} \dots B^{j_{r}}_{i_{2r-1}i_{2r}}.$$
\end{theorem}
\noindent Again, all of this is in \cite{GM}.

For an example, let $V$ have dimension two and $W$ have dimension 
one.  Let $r=1$.  Then our vector-valued form $B$ can be 
represented either as a 
two form
$$ a e^{1}\otimes e^{1}+ b  e^{1}\otimes e^{2} + c  e^{2}\otimes e^{1} 
+ d e^{2}\otimes e^{2}$$
or as a two by two matrix
$$\left(
\begin{array}{cc}
 a & b \\
 c & d
\end{array} \right).$$
We  set 
$$v_{\sigma} = \varepsilon^{I} (2,2,\mbox{identity} )e_{i_{1}}
 \otimes e_{i_{2}}
 =  e_{1}\otimes e_{2} - e_{2}\otimes e_{1} .$$
 Since $W$ has dimension one, we must have $w^{\eta}$ be the identity.
 Then $v_{\sigma}\otimes w^{\eta}$ acting on $B$ will be
 $$b - c$$ and is zero precisely when  the matrix $B$ 
 is symmetric.

\section{Relations among invariants for the general linear group}

\subsection{Nontrivial relations}

For $Gl(n,\mbox{\bf {C}})$ acting on a vector space $V$, classically  not only are the invariants known,
 but also  so are the relations. 
 Everything in this section is well-known. We will first discuss the
 second fundamental theorem in the language
of brackets, or determinants.  This is the invariant language for which the 
second fundamental theorem is the most clear.  We then will state the second fundamental
theorem for the two cases that we need in this paper, namely for $V\otimes V$ and
 $W^{*}$.

Let ${\bf v}_1, \ldots, {\bf v}_n$ be $n$ column vectors in ${\bf C}^n$.  The 
general linear group $Gl(n,{\bf C})$ acts on the vectors in  ${\bf C}^n$ by multiplication
on the left.  Classically, the {\it bracket} of the vectors ${\bf v}_1, \ldots, {\bf v}_n$,
denoted by $[{\bf v}_1, \ldots, {\bf v}_n]$, is defined to be the determinant of 
the $n\times n$ matrix whose columns are the vectors ${\bf v}_1, \ldots, {\bf v}_n$.
Thus by definition
$$ [{\bf v}_1, \ldots, {\bf v}_n] = \det ({\bf v}_1, \ldots, {\bf v}_n).$$
By basic properties of the determinant, we have that the bracket is an invariant, since
\begin{eqnarray*} [A{\bf v}_1, \ldots, A{\bf v}_n]&=& \det (A{\bf v}_1, \ldots,A{\bf v}_n)\\
&=& \det(A)\det ({\bf v}_1, \ldots, {\bf v}_n) \\
&=& \det(A)[{\bf v}_1, \ldots, {\bf v}_n]
\end{eqnarray*}
The punchline of the first fundamental theorem in this language is that the only 
invariants are combinations of various brackets and hence of determinants. (See p. 22 in 
\cite{DC} or p. 45 in \cite{Weyl}.)

The second fundamental theorem  reflects the fact that the determinant of a matrix with
two identical rows is zero. Choose $n+1$ column vectors ${\bf v}_1, \ldots, {\bf v}_{n+1}$
 in ${\bf C}^n$. Label the entries of the vector ${\bf v}_i=(v_{ij})$, for $1\leq j \leq n$.
Let ${\bf e}_k$ denote the vectors in the standard basis for  ${\bf C}^n$.  Thus all of 
the entries in  ${\bf e}_k$ are zero, except in the $k$th entry, which is one.
Then 
\begin{eqnarray*}[{\bf v}_i, {\bf e}_2, \ldots, {\bf e}_n]&=& \det ({\bf v}_i, {\bf e}_2, \ldots, {\bf e}_n)\\
&=& v_{i1}
\end{eqnarray*}
Now consider the $(n+1)\times (n+1)$ matrix
$$V= \left(
\begin{array}{cccc}
 v_{11} & v_{21}& \cdots & v_{(n+1)1} \\
 v_{11} & v_{21}& \cdots & v_{(n+1)1} \\
 v_{12} & v_{22}& \cdots & v_{(n+1)2}  \\
\vdots & \vdots & \vdots & \vdots \\
 v_{1n} & v_{2n}& \cdots & v_{(n+1)n}
\end{array} \right).$$
Since its top two rows are identical, its determinant is zero.  Then
\begin{eqnarray*}0 &=& \det(V)\\
&=& \sum_{k=1}^{n+1} 
(-1)^{k+1} v_{k1}\det({\bf v}_1, \ldots, \hat{{\bf v}}_k, \ldots, {\bf v}_{n+1})\\
&=&  \sum_{k=1}^{n+1} 
(-1)^{k+1} [{\bf v}_1, \ldots, \hat{{\bf v}}_k, \ldots, {\bf v}_{n+1}]
[{\bf v}_k, {\bf e}_2, \ldots, {\bf e}_n],
\end{eqnarray*}
where $ \hat{{\bf v}}_k$ means delete the $ {\bf v}_k$ term.  This equation is a relation
among brackets.  The punchline of the second fundamental theorem is that all nontrivial 
relations
are of the form
$$\sum_{k=1}^{n+1} 
(-1)^{k+1} [{\bf v}_1, \ldots, \hat{{\bf v}}_k, \ldots, {\bf v}_{n+1}]
[{\bf v}_k, {\bf w}_2, \ldots, {\bf w}_n] = 0, $$
where ${\bf w}_2, \ldots, {\bf w}_{n}$ are any column
vectors. 

 We will  define {\it nontrivial} in the next subsection.  Basically 
 the trivial relations  stem from the fact that  rearranging the columns of a matrix will only
change the determinant by at most a sign.  Thus rearranging the vectors
${\bf v}_1, \ldots, {\bf v}_{n}$ and then taking the bracket
will give us, up to sign, the same invariant.

Now to  quickly put this into the language of tensors, first for 
the relations for $V\otimes V$.  Using the 
notation of the last section, let $n$ denote the dimension of $V$ and let $2r=nu$.  We know that all invariants 
for a given $r$
are generated by all possible $v_{\sigma}
 = \varepsilon^{I}(n, 2r, \sigma) e_{i_{1}}\otimes \ldots 
\otimes e_{i_{2r}}$.  All of these invariants have the same weight.   
Thus, for each $r$, the sum of any two $v_{\sigma} $  spans another 
invariant one dimensional subspace in  $(V\otimes V)^{\otimes r}$.  Denote the subspace 
of  $(V\otimes V)^{\otimes r}$ spanned by the varous $v_{\sigma}$, 
with $\sigma \in S_{2r}$, by
$$(V\otimes V)^{r}_{\mbox inv}.$$
Then the first fundamental theorem in this context can be interpreted 
as giving a spanning set for $(V\otimes V)^{r}_{\mbox inv}$, for 
each $r$.   For each $r$, denote the vector space with basis indexed by
the $v_{\sigma}$ by $$(V\otimes V)^{r}_{0}.$$  There is thus an onto 
linear transformation 
$$(V\otimes V)^{r}_{0} \rightarrow (V\otimes V)^{r}_{\mbox inv}.$$
We want to find a spanning set of the kernel of this map.  This will 
be a set of relations among the generators of the 
invariants.

We first need some notation.  Let $\{i_{1},\ldots, 
i_{n+1}\} $ be a subset of $n+1$ distinct elements chosen from 
$\{1,2, \ldots, 2r\}$ and let $\sigma \in S_{2r}$. For $1\leq j\leq n+1$, define
 $\sigma_j \in S_{2r}$ by setting
$\sigma_j (k) = \sigma(k)$ if $k$ is not in $ \{i_{1},\ldots, 
i_{n+1}\} $ and
$$\sigma_j (i_k) = \left\{ \begin{array}{ll}
				\sigma (i_k) & \mbox{if}\; k<j\\
				\sigma(i_{k-1}) & \mbox{if}\; j<k \\
				\sigma(i_{n+1})& \mbox{if}\; j=k
\end{array}\right. $$
For a fixed $\sigma$ and subsequence $\{i_{1},\ldots, 
i_{n+1}\}$, let $A(\sigma,\{i_{1},\ldots, 
i_{n+1}\}) $ denote the subset of  $S_{2r}$ consisting of the $\sigma_j$.
Then the earlier stated second fundamental theorem can be reformulated 
in this context as:

\begin{theorem}
All nontrivial relations for $(V\otimes V)^{\otimes r}$ 
are linear combinations of 
$$\sum_{\sigma_j \in
A(\sigma,\{i_{1},\ldots, 
i_{n+1}\})}(-1)^{j+1}v_{\sigma_j}=0,$$
for all possible subsets $\{i_{1},\ldots, 
i_{n+1}\}$ and all possible $\sigma \in S_{2r}$.
\end{theorem}
\noindent  Again, the term {\it nontrivial} means the same as before and is 
hence 
 simply dealing with the fact that if 
you flip two columns of a matrix, the corresponding determinants 
changes sign. 

Note that this theorem states that the relations 
are linear for generators $v_{\sigma}$ with $\sigma \in S_{2r}$.  
Thus we are indeed capturing a spanning set for the kernel of 
the map $(V\otimes V)^{r}_{0} \rightarrow (V\otimes V)^{r}_{\mbox 
inv},$ for each r.   
Fixing $r$, denote the vector space with basis indexed 
by each of the above relations and by each of the trivial relations by $$(V\otimes V)^{r}_{1}.$$  Then we 
have an exact sequence
$$(V\otimes V)^{r}_{1} \rightarrow(V\otimes V)^{r}_{0} \rightarrow (V\otimes 
V)^{r}_{\mbox 
inv},$$
an exact sequence that is a linear algebra description of both the 
first and second fundamental theorems for this particular group 
action.

The relations for the invariants of the general linear group $Gl(k,C)$ 
acting on $W^{*}$ are similar.  Here we have $r=kv$.  The invariants are generated by 
$w^{\eta}= \varepsilon_{J}(k, r, \eta) f^{i_{1}}\otimes \ldots 
\otimes f^{i_{r}}.$ Choose $k+1$ distinct elements $\{i_{1},\ldots, 
i_{k+1}\} $  from $\{1,2, \ldots, r\}$ and an $\eta \in S_r$.  Let 
$B(\eta, \{i_{1},\ldots, 
i_{k+1}\})$ denote the set of all  $\eta_j \in S_{r}$, for $1\leq j\leq k+1$, defined by setting
$\eta_j (l) = \eta(l)$ if $i$ is not in $ \{i_{1},\ldots, 
i_{k+1}\} $ and
$$\eta_j (i_l) = \left\{ \begin{array}{ll}
				\eta (i_l) & \mbox{if}\; l<j\\
				\eta(i_{l-1}) & \mbox{if}\; j<l \\
				\eta(i_{n+1})& \mbox{if}\; j=l
\end{array}\right. $$
Then
\begin{theorem}
All nontrivial relations for $(W^{*})^{\otimes r}$ 
are linear combinations of 
$$\sum_{\eta_j \in
B(\eta,\{i_{1},\ldots, 
i_{k+1}\}) } (-1)^jw^{\eta_j}=0,$$
for all possible  $\{i_{1},\ldots, 
i_{k+1}\} $  and  $\eta \in S_r$.
\end{theorem}
The proofs of theorems 3 and 4 are in \cite{Weyl} on pp. 70-76.  Weyl uses the bracket notation, but
the equivalence is straightforward.  A matrix approach is in section II.3, on page 71, in
\cite{ACGH}.

Mirroring what we did above, we know that the invariant linear 
subspaces  for $(W^{*})^{\otimes r}$ are generated by all possible 
$w^{\eta}$, for $\eta \in S_r$ and that all of these invariants 
have the same weight.  Thus for each $r$, the various 
$w^{\eta}$ span an invariant subspace of $(W^{*})^{\otimes r}$.  Denote 
this subspace by 
$$(W^{r})_{\mbox inv}.$$
 For each $r$, let the vector space with basis indexed by
the $w^{\eta}$  be  $(W^{*})^{r}_{0}.$
Let the vector space with basis indexed by the above relations for the 
various $w^{\eta}$ and by the trivial relations be denoted by $(W^{*})^{r}_{1}.$  Then we have the exact 
sequence
$$(W^{*})^{r}_{1} \rightarrow (W^{*})^{r}_{0} \rightarrow (W^{*})^{r}_{\mbox 
inv},$$

Now for an example.  We first will write down a relation in the bracket notation, give the 
translation in terms of tensors and then see that this explicit relation is in the above list.
Let ${\bf v}_1,{\bf v}_2, {\bf v}_3, {\bf w}$ be any four column vectors in {\bf C}$^2$.
Then by explicit calculation we have
$$[{\bf v}_1, {\bf v}_2][{\bf v}_3,{\bf w}]-[{\bf v}_1, {\bf v}_3][{\bf v}_2,{\bf w}] +
[{\bf v}_2, {\bf v}_3][{\bf v}_1,{\bf w}]=0.$$
In the dimension two vector space $W^{*}$, with basis $f^1,f^2$, consider the corresponding relation
\begin{eqnarray*}
\Sigma&=& w^{(1)}-w^{(23)}+ w^{(132)} \\
&=& (f^{1}\otimes f^{2}\otimes f^{1}\otimes f^{2} -f^{1}\otimes f^{2}\otimes f^{2}\otimes f^{1}\\
  &&   -f^{2}\otimes f^{1}\otimes f^{1}\otimes f^{2} +f^{2}\otimes f^{1}\otimes f^{2}\otimes
     f^{1}) \\
&& +(f^{2}\otimes f^{1}\otimes f^{1}\otimes f^{2}
     -f^{2}\otimes f^{2}\otimes f^{1}\otimes f^{1}\\
  &&   -f^{1}\otimes f^{1}\otimes f^{2}\otimes f^{2} +f^{1}\otimes f^{2}\otimes f^{2}\otimes
     f^{1}) \\
&&+(f^{1}\otimes f^{1}\otimes f^{2}\otimes f^{2}
     -f^{2}\otimes f^{1}\otimes f^{2}\otimes f^{1}\\
&&     -f^{1}\otimes f^{2}\otimes f^{1}\otimes f^{2} +f^{2}\otimes f^{2}\otimes f^{1}\otimes
     f^{1}) \\
 & = & 0.
\end{eqnarray*}
Now to show that the relation $\Sigma$ is in the above list.
We have $r=4$.  Let $\eta \in S_4$ be the identity permutation.  Let our sequence 
$\{i_1,i_2, i_3\}$ be simply $\{1,2,3\}$.  Then $\eta_1$ is the permutation $(132)$, $\eta_2$
is the permutation $(23)$ and $\eta_3$ is the identity permutation.  Thus the relation $\Sigma$
is an example of the relation:
$$ w^{\eta_3}-w^{\eta_2}+ w^{\eta_1}=0.$$

\subsection{Trivial relations}

All of this section is still classical.  

It can be directly checked, continuing with our example for the two 
dimensional vector space $W^{*}$, that
$$w^{(1)}+w^{(123)}+ w^{(132)}=0.$$
Here the invariant $w^{(123)}$ is trivially related to the invariant $w^{(23)}$ (specifically,
$w^{(123)}=-w^{(23)}$).  This is easiest to see in bracket notation, as this is just reflecting that
$$[{\bf v}_3, {\bf v}_1][{\bf v}_2,{\bf w}]=-[{\bf v}_1, {\bf v}_3][{\bf v}_2,{\bf w}],$$
which in turn simply reflects that fact that the sign of a determinant changes when we flip two
columns.

This is the source of all relations that we want to call {\it 
trivial}.   Rearranging the columns of a matrix will not change the determinant 
if the rearrangement is given by an even permutation of the 
permutation group and will change the determinant by a sign
if the rearrangement is given by an odd permutation of the 
permutation group.    

 We will give the explicit criterion for trivial relations for the 
 case of a $k$ dimensional vector space $W^{*}$. As always, let $r= 
 kv$. Our 
 goal is to determine, given $\sigma, \tau \in S_{r}$, when 
 $$w^{\sigma}=\pm w^{\tau}.$$
 This happens when we have the following $v$ equalities of sets:
 \begin{eqnarray*}
 \{ \sigma^{-1}(1),\ldots, \sigma^{-1}(k)\}&=& \{ \tau^{-1}(1),\ldots, 
 \tau^{-1}(k)\}\\
   \{ \sigma^{-1}(k+1),\ldots, \sigma^{-1}(2k)\}&=& \{ \tau^{-1}(k+1),\ldots, 
 \tau^{-1}(2k)\}\\
 &\vdots& \\
  \{ \sigma^{-1}((v-1)k+1),\ldots, \sigma^{-1}(kv)\}&=& \{ \tau^{-1}((v-1)k+1),\ldots, 
 \tau^{-1}(kv)\}
 \end{eqnarray*}
 Each set on the right is thus a permutation of the corresponding set 
 on the left. 
 We will have $w^{\sigma}= w^{\tau}$ if there are an even number of 
 odd permutations taking the left hand side of the above set 
 equalities to the right and $w^{\sigma}= -w^{\tau}$  if there are an 
 odd number.
 
 Consider our initial example $w^{(123)}=-w^{(23)}$ when $W^{*}$ is 
 two dimensional.  Let $\sigma = (123)$ and $\tau=(23)$.
 Then
 $$\sigma^{-1}(1)= 3, \sigma^{-1}(2)= 1, \sigma^{-1}(3)= 2,
 \sigma^{-1}(4)= 4$$
 and
 $$\tau^{-1}(1)= 1, \tau^{-1}(2)= 3, \sigma^{-1}(3)= 2, 
 \sigma^{-1}(4)= 4.$$
 Then $\{\sigma^{-1}(1), \sigma^{-1}(2)\}$ is an odd permutation of 
 $\{\tau^{-1}(1), \tau^{-1}(2)\}$, while $\{\sigma^{-1}(3), 
 \sigma^{-1}(4)\}$ is exactly the same as $\{\tau^{-1}(1), 
 \tau^{-1}(2)\}$, reflecting the fact that $w^{(123)}=-w^{(23)}$.

 \section{A Second Fundamental Theorem for \\
 $(V\otimes V\otimes W^{*})^{\otimes r}$}
 We have the two exact sequences
 
 $$(V\otimes V)^{r}_{1} \rightarrow(V\otimes V)^{r}_{0} \rightarrow (V\otimes 
V)^{r}_{\mbox 
inv}$$
 and
 $$(W^{*})^{r}_{1} \rightarrow (W^{*})^{r}_{0} \rightarrow ((W^{*})^{r})_{\mbox{inv}}.$$
 Tensoring either of these exact sequences by a complex vector space 
 will maintain the exactness.
The point of \cite{GM} is that the natural map
$$(V\otimes V)^{r}_{0}\otimes  (W^{*})^{r}_{0} \rightarrow (V\otimes 
V)^{r}_{\mbox{inv}}\otimes (W^{*})^{r}_{\mbox{inv}}$$
is onto.  We want to find the kernel of this map.

We  have the following commutative double exact sequece:
$$\begin{array}{ccccccc}
0        & & 0        & & 0        & & \\
\uparrow & & \uparrow & & \uparrow & &  \\
(V\otimes V)^{r}_{1} \otimes (W^{*})^{r}_{\mbox{inv}}& \rightarrow & (V\otimes V)^{r}_{0}
 \otimes (W^{*})^{r}_{\mbox{inv}}& \rightarrow & (V\otimes 
 V)^{r}_{\mbox{inv}} \otimes (W^{*})^{r}_{\mbox{inv}}& \rightarrow& 0\\
 \uparrow & & \uparrow & & \uparrow & &  \\
 (V\otimes V)^{r}_{1} \otimes (W^{*})^{r}_{0}& \rightarrow & (V\otimes V)^{r}_{0}
 \otimes (W^{*})^{r}_{0}& \rightarrow & (V\otimes 
 V)^{r}_{\mbox{inv}} \otimes (W^{*})^{r}_{0}& \rightarrow& 0\\
 \uparrow & & \uparrow & & \uparrow & &  \\
 (V\otimes V)^{r}_{1} \otimes (W^{*})^{r}_{1}& \rightarrow & (V\otimes V)^{r}_{0}
 \otimes (W^{*})^{r}_{1}& \rightarrow & (V\otimes 
 V)^{r}_{\mbox{inv}} \otimes (W^{*})^{r}_{1}& \rightarrow& 0\\
 \end{array}$$

 \noindent with 
 $$(V\otimes V)^{r}_{0}
 \otimes (W^{*})^{r}_{0} \rightarrow (V\otimes 
 V)^{r}_{\mbox{inv}} \otimes (W^{*})^{r}_{\mbox{inv}}$$
 from the above double exact sequence being onto.  All of the above 
 maps are linear transformations of vector spaces.  A second fundamental 
 theorem for vector-valued forms will be a description of the kernel 
 of this map.
 
\begin{theorem}
Under the natural maps from the above double exact sequence, the kernal of the map from $(V\otimes V)^{r}_{0}
 \otimes (W^{*})^{r}_{0}$
 to $(V\otimes 
 V)^{r}_{\mbox{inv}} \otimes (W^{*})^{r}_{\mbox{inv}}$ is 
 $$(V\otimes V)^{r}_{1}
 \otimes (W^{*})^{r}_{0} \oplus (V\otimes V)^{r}_{0}
 \otimes (W^{*})^{r}_{1}.$$
 \end{theorem}
The proof is a routine diagram chase.
 
Thus by standard arguments involving commutative diagrams,   
the following 
 sequence of vector spaces is exact:

 \begin{eqnarray*} (V\otimes V)^{r}_{1}
 \otimes (W^{*})^{r}_{0} \oplus (V\otimes V)^{r}_{0}
 \otimes (W^{*})^{r}_{1} &\rightarrow& (V\otimes V)^{r}_{0}
 \otimes (W^{*})^{r}_{0}\\
 & \rightarrow& (V\otimes 
 V)^{r}_{\mbox{inv}} \otimes (W^{*})^{r}_{\mbox{inv}} \\
 &\rightarrow& 0.
 \end{eqnarray*}

In another language, this theorem can be stated as:

 \begin{theorem}[A Second Fundamental Theorem] Among invariants of vector-valued bilinear forms,
there exist relations of the following type:  $$\left( (V
\otimes V)_{1} \otimes (W^{*})_{0} \right) \oplus \left( (V
\otimes V)_{0} \otimes (W^{*})_{1} \right) .$$
All relations are linear combinations of the above relations.
\end{theorem}

\section{An example of a relation}
Let $V$ and $W^{*}$ both have dimension two.  Recall our example of a relation for
$W$:
$$\Sigma = w^{1} -w^{(23)} + w^{(132)}.$$
Here $k=2$ and $r=4$.  Choose $(23)(67)\in S_8$.  Then

$$v_{(23)(67)}  = \varepsilon^{I}(2, 8, (23)(67)) e_{i_{1}}\otimes \ldots 
\otimes e_{i_{2r}}$$
\begin{eqnarray*}&=&
e_1 \otimes e_1 \otimes e_2 \otimes e_2 \otimes e_1 \otimes e_1 \otimes e_2 \otimes e_2 
- e_1 \otimes e_1 \otimes e_2 \otimes e_2 \otimes e_1 \otimes e_2 \otimes e_2 \otimes e_1\\
&&-e_1 \otimes e_1 \otimes e_2 \otimes e_2 \otimes e_2 \otimes e_1 \otimes e_1 \otimes e_2
+e_1 \otimes e_1 \otimes e_2 \otimes e_2 \otimes e_2 \otimes e_2 \otimes e_1 \otimes e_1\\
&&-e_1 \otimes e_2 \otimes e_2 \otimes e_1 \otimes e_1 \otimes e_1 \otimes e_2 \otimes e_2
+e_1 \otimes e_2 \otimes e_2 \otimes e_1 \otimes e_1 \otimes e_2 \otimes e_2 \otimes e_1\\
&&+ e_1 \otimes e_2 \otimes e_2 \otimes e_1 \otimes e_2 \otimes e_1 \otimes e_1 \otimes e_2
-e_1 \otimes e_2 \otimes e_2 \otimes e_1 \otimes e_2 \otimes e_2 \otimes e_1 \otimes e_1\\
&&-e_2 \otimes e_1 \otimes e_1 \otimes e_2 \otimes e_1 \otimes e_1 \otimes e_2 \otimes e_2
+e_2 \otimes e_1 \otimes e_1 \otimes e_2 \otimes e_1 \otimes e_2 \otimes e_2 \otimes e_1\\
&&+e_2 \otimes e_1 \otimes e_1 \otimes e_2 \otimes e_2 \otimes e_1 \otimes e_1 \otimes e_2
-e_2 \otimes e_1 \otimes e_1 \otimes e_2 \otimes e_2 \otimes e_2 \otimes e_1 \otimes e_1\\
&&+e_2 \otimes e_2 \otimes e_1 \otimes e_1 \otimes e_1 \otimes e_1 \otimes e_2 \otimes e_2
-e_2 \otimes e_2 \otimes e_1 \otimes e_1 \otimes e_1 \otimes e_2 \otimes e_2 \otimes e_1\\
&&-e_2 \otimes e_2 \otimes e_1 \otimes e_1 \otimes e_2 \otimes e_1 \otimes e_1 \otimes e_2
+e_2 \otimes e_2 \otimes e_1 \otimes e_1 \otimes e_2 \otimes e_2 \otimes e_1 \otimes e_1.
\end{eqnarray*}

Then  we have the relation
\begin{eqnarray*}v_{(23)(67)}\otimes \Sigma &= &v_{(23)(67)}\otimes w^{(1)}
-v_{(23)(67)}\otimes w^{(23)} + v_{(23)(67)}\otimes w^{(132)}\\
&=& 0,
\end{eqnarray*}
which can now be directly checked.

\section{Finding Relations for Bilinear Forms}

Most people, though, are not that interested in invariant one-dimensional
subspaces of  $(V\otimes V\otimes W^{*})^{\otimes r}$ for various
positive integers $r$, 
 but are more interested in invariants of bilinear forms.  As 
 discussed in section 1.3, this means that we are interested in the 
 algebra  of homogeneous  polynoials in 
$\mathbb{C}[V^{*}\otimes V^{*}\otimes W]$ that are invariant under 
the previously defined group action by $Aut(V) \times Aut(W)$.  But 
the homogoenous polynomials of degree $r$
 can be identified to elements in the symmetric space 
 $(V^{*}\otimes V^{*}\otimes W)^{\odot  r}$.  We thus want to find 
 the invariant lines in the dual space and hence the invariant one-dimensional subspaces of 
$(V\otimes V\otimes W^{*})^{\odot  r}$.  So far all we have are the invariant one dimensional
subspaces of  $(V\otimes V\otimes W^{*})^{\otimes r}$.

Denote the vector space spanned by the invariant one-dimensional subspaces 
in $(V\otimes V\otimes W^{*})^{\odot  r}$  by
$$(V\otimes V\otimes W^{*})^{\odot  r}_{\mbox{inv}}.$$
There is  a natural onto map from $(V\otimes 
V)^{r}_{\mbox{inv}}
 \otimes W^{r}_{\mbox{inv}}$ to $(V^{*}\otimes V^{*}\otimes W)^{\odot  
 r}_{\mbox{inv}}$.  This is simply  the restriction to $(V\otimes 
V)^{r}_{\mbox{inv}}
 \otimes W^{r}_{\mbox{inv}}$ of the 
 symmetrizing map 
 $$S:(V\otimes V\otimes W^{*})^{\otimes r}\rightarrow (V\otimes V\otimes 
 W^{*})^{\odot r}.$$
 We need, though,  to check that an element of $(V\otimes V\otimes W^{*})^{\otimes 
 r}$ that generates a one-dimensional  invariant subspace still 
 generates a one-dimensional invariant subspace after the application 
 of the map $S$.

 The permutation group $S_{r}$ acts naturally on both $(W^{*})^{\otimes r}$ 
 and on $(V\otimes V)^{\otimes r}$.  Let $\tau\in S_{r}$.  Then for 
 any $\eta \in S_{r}$, it can be directly checked that 
 $$\tau(w^{\eta})=w^{\eta\cdot \tau^{-1}},$$
 another element in our list of generators.
 
 Similarly, for any $\tau\in S_{r}$, we will have $\tau(v_{\sigma})$ be in our list of 
 generators,  for any $v_{\sigma}$.  Here the notation is a bit 
 cumbersome.  Each $\tau\in S_{r}$ will induce an element 
 $\hat{\tau}\in S_{2r}$, where if $\tau(i)=j$, then 
 \begin{eqnarray*}\hat{\tau}(2i-1)&=&2j-1\\ 
 \hat{\tau}(2i)&=&2j.
 \end{eqnarray*}
   Then it can be 
 checked that 
 $$\tau(v_{\sigma})=v_{\sigma \cdot \hat{\tau}^{-1}}.$$
 Thus $\tau(v_{\sigma})$ is another invariant.

Then,  as is implicit in \cite{GM}, we have:
\begin{theorem}[First Fundamental Theorem for  Bilinear Forms] 
The vector space $ (V\otimes V \otimes W^{*})^{\odot r} $ 
has an invariant one-dimensional subspace if and 
only if $n$ divides $2r$ and $k$ 
divides $r$.  Every invariant one-dimensional subspace is a linear 
combination of various $S(v_{\sigma}\otimes w^{\eta})$, where $\sigma$ 
and $\eta$ range   through $S_{2r}$ and $S_{r}$, respectively.  
\end{theorem}

We are interested in the relations among the various $S(v_{\sigma}\otimes w^{\eta})$.
Since $(V\otimes V\otimes 
 W^{*})^{\odot r}$ is contained in $(V\otimes V\otimes W^{*})^{\otimes r}$, 
 any relation in $(V\otimes V\otimes 
 W^{*})^{\odot r}$  must be in the relations for $(V\otimes V\otimes 
 W^{*})^{\odot r}$.

Hence  we just need to map all of our previous relations 
 to $(V\otimes V\otimes 
 W^{*})^{\odot r}$ via $S$.   If there is a relation of invariants in $(V\otimes V\otimes 
 W^{*})^{\odot r}$, we have already captured it.

 Thus we have
 
 \begin{theorem}[Second Fundamental Theorem for  Bilinear Forms]
 All nontrivial relations among the nonzero elements  $S(v_{\sigma}\otimes 
 w^{\eta})$, where $\sigma$ 
and $\eta$ range   through $S_{2r}$ and $S_{r}$, respectively, are 
linear combinations of all
$$S(v_{\sigma}\otimes \sum_{\eta_j \in
B(\eta,\{i_{1},\ldots, 
i_{k+1}\}) } (-1)^jw^{\eta_j})$$
and 
$$S(\sum_{\sigma_j \in
A(\sigma,\{i_{1},\ldots, 
i_{n+1}\})}(-1)^{j+1}v_{\sigma_j}\otimes w^{\eta}).$$
  \end{theorem}

We had to use the term ``nonzero'' in the above theorem.    Some of our invariants in $(V\otimes V\otimes 
 W^{*})^{\otimes r}$ will be mapped to zero under $S$.  None of the above 
 describes the      kernel of $S$. As we will see in the next 
 section, this does happen.  In fact, if we consider the example of 
 symmetrizing
 map $S:(W^{*})^{\otimes r}\rightarrow (W^{*})^{\odot r},$ then it 
 is not at all obvious that every $v_{\sigma}\otimes w^{\eta}$ is not 
 sent to zero, since 
 $$S(w^{\eta})=0$$
 for all $\eta \in S_{r}$, which can be directly checked.  (This just reflects 
 that the geometric fact that there 
 are no invariants for a singe vector in a vector space under the 
 group action of the automorphisms of the vector space, since any 
 vector can be sent to any other vector.)  Again, we will see an example 
 in the next section that there are nontrivial relations.
 
 Thus we have an algorithm for finding the invariants and for finding 
 the relations. For each r, we just map all of our generators in $(V\otimes V\otimes 
 W^{*})^{\otimes r}$ to $(V\otimes V\otimes 
 W^{*})^{\odot r}$ by $S$, disposing of those that map to zero.  All the
 remaining relations will be already be accounted for by applying $S$ 
 to the previous relations.

 \section{An Example}
 
  We now translate  the above relations into the
language of invariant polynomials of bilinear forms, for the 
particular case of an element in  $(V \otimes V \otimes W^{*})^{\otimes 4}$. Let
$$B= \left(\left(
\begin{array}{cc}
 b_{11}^1 & b_{12}^1 \\
 b_{21}^1 & b_{22}^1
\end{array} \right),\left(\begin{array}{cc}
 b_{11}^2 & b_{12}^2 \\
 b_{21}^2 & b_{22}^2
\end{array} \right)\right)$$
be a bilinear form from $V\times V$ to $W$.  As a tensor in $V^*\otimes V^* \otimes W$, this 
bilinear form becomes the tensor:
\begin{eqnarray*}b&=& b_{ij}^k e^i\otimes e^j  \otimes f_k\\
&=&  b_{11}^1 e^1\otimes e^1 \otimes f_1 +  b_{12}^1 e^1\otimes e^2 \otimes f_1
+ b_{21}^1 e^2\otimes e^1 \otimes f_1 +  b_{22}^1 e^2\otimes e^2 \otimes f_1\\
&&+ b_{11}^2 e^1\otimes e^1 \otimes f_2 +  b_{12}^2 e^1\otimes e^2 \otimes f_2
+ b_{21}^2 e^2\otimes e^1 \otimes f_2 +  b_{22}^1 e^2\otimes e^2 \otimes f_2.
\end{eqnarray*}
To see what the invariants $v_{(23)(67)}\otimes w^{(1)}$,
$v_{(23)(67)}\otimes w^{(23)}$ and $v_{(23)(67)}\otimes w^{(132)}$  
are in terms of the variables
$b_{ij}^k$, we have each act on the tensor $b\otimes b\otimes  b \otimes b$.
Then, after a painful calculation, we get that
\begin{eqnarray*}
v_{(23)(67)}\otimes w^{(1)}(b\otimes b\otimes  b \otimes b)&=&0\\
v_{(23)(67)}\otimes w^{(23)}(b\otimes b\otimes  b \otimes b)&=&
4b_{11}^1 b_{22}^1 b_{11}^2 b_{22}^2 + 4 b_{12}^1 b_{22}^1 b_{11}^2 b_{21}^2\\
&&-8 b_{12}^1 b_{21}^1 b_{11}^2 b_{22}^2 -8 b_{11}^1 b_{22}^1 b_{12}^2 b_{21}^2\\&&
+ 4 b_{11}^1 b_{21}^1 b_{12}^2 b_{22}^2 + 4 b_{21}^1 b_{22}^1 b_{11}^2 b_{12}^2\\&&
-2 b_{22}^1 b_{22}^1 b_{11}^2 b_{11}^2 -2 b_{11}^1 b_{11}^1 b_{22}^2 b_{22}^2\\&&
+4 b_{12}^1 b_{21}^1 b_{12}^2 b_{21}^2 -2 b_{21}^1 b_{21}^1 b_{12}^2 b_{12}^2\\&&
-2 b_{12}^1 b_{12}^1 b_{21}^2 b_{21}^2 +4 b_{11}^1 b_{12}^1 b_{21}^2 b_{22}^2\\
v_{(23)(67)}\otimes w^{(132)}(b\otimes b\otimes  b \otimes b)&=&
4b_{11}^1 b_{22}^1 b_{11}^2 b_{22}^2 + 4 b_{12}^1 b_{22}^1 b_{11}^2 b_{21}^2\\&&
-8 b_{12}^1 b_{21}^1 b_{11}^2 b_{22}^2 -8 b_{11}^1 b_{22}^1 b_{12}^2 b_{21}^2\\&&
+ 4 b_{11}^1 b_{21}^1 b_{12}^2 b_{22}^2 + 4 b_{21}^1 b_{22}^1 b_{11}^2 b_{12}^2\\&&
-2 b_{22}^1 b_{22}^1 b_{11}^2 b_{11}^2 -2 b_{11}^1 b_{11}^1 b_{22}^2 b_{22}^2\\&&
+4 b_{12}^1 b_{21}^1 b_{12}^2 b_{21}^2 -2 b_{21}^1 b_{21}^1 b_{12}^2 b_{12}^2\\&&
-2 b_{12}^1 b_{12}^1 b_{21}^2 b_{21}^2 +4 b_{11}^1 b_{12}^1 b_{21}^2 b_{22}^2
\end{eqnarray*}
Since $v_{(23)(67)}\otimes w^{(23)}(b\otimes b\otimes  b \otimes b)=
v_{(23)(67)}\otimes w^{(132)}(b\otimes b\otimes  b \otimes b)$ and since
$v_{(23)(67)}\otimes w^{(1)}(b\otimes b\otimes  b \otimes b)=0$, we see that our relations
from the tensor language do translate to relations on the invariant polynomials in the 
$b_{ij}^k$ terms.

Also, note that this means that 
$$S(v_{(23)(67)}\otimes w^{(1)})=0$$
and hence $v_{(23)(67)}\otimes w^{(1)}$ is in the kernel of the 
symmetrizing map $S:(V\otimes V\otimes W^{*})^{\otimes r}\rightarrow 
(V\otimes V\otimes W^{*})^{\odot r}$.

Of course, it is difficult to see what these invariants and relations 
actually measure.  There must be geometry behind these formulas, 
though it is almost always hidden.  There are times that we can 
understand some of the information contained in the formulas.  As an 
example, we 
will now  see why we chose $r=4$ for our example.  We will see that for when 
the rank of $V$ and 
$W$ are both two, this is the first time we would expect any interesting invariants
for bilinear forms.
 Consider the relation 
$$v_{(1)} -v_{(23)} +v_{(132)}=0$$
for  $V^{*}\otimes 
V^{*}$. (This is the $V$-analogue of our earlier $w^{1} -w^{(23)} + w^{(132)}=0$.)
Next we construct an invariant  of $(W)_0$ by letting $k=2$, $r=2$,
and $\eta = \mbox{id}$.  This gives us the invariant $w^{\eta}= f^{1}\otimes
f^{2} -f^{2} \otimes f^{1}$. Then a relation is 
$$(v_{(1)} -v_{(23)} +v_{(132)})\otimes w^{\eta} =0.$$
But this relation is not at all interesting when made into a statement about polynomials on the
entries of a bilinear form, since each
invariant becomes the zero polynomial:
\begin{eqnarray*}(v_{(1)}\otimes  w^{\eta})(b\otimes b)&=&0\\
(v_{(23)}\otimes  w^{\eta})(b\otimes b)&=&0\\
(v_{(132)}\otimes  w^{\eta})(b\otimes b)&=&0,
\end{eqnarray*}
following from a direct calculation.

 When $n=k=2$, the first time that interesting invariants can occur is indeed when $r=4$ (which is why
our example has $r=4$), as we will now see. 
 Start with our bilinear form
$$B= (B_1,B_2)= \left(\left(
\begin{array}{cc}
 b_{11}^1 & b_{12}^1 \\
 b_{21}^1 & b_{22}^1
\end{array} \right),\left(\begin{array}{cc}
 b_{11}^2 & b_{12}^2 \\
 b_{21}^2 & b_{22}^2
\end{array} \right)\right)$$
and consider the polynomial
\begin{eqnarray*} P(x,y)& = &\det(xB_1 + yB_2) \\
&=& \det \left(\begin{array}{cc}
x b_{11}^1 + y b_{11}^2 &x b_{12}^1 + y b_{12}^2\\
 xb_{21}^1 + y b_{21}^2 &x b_{22}^1 +y b_{22}^2\\
\end{array} \right)\\
&=&(b_{11}^{1}b_{22}^{1}-b_{12}^1 b_{21}^1)x^2 +(b_{11}^1 b_{22}^2 + b_{22}^1 b_{11}^2
-b_{12}^1 b_{21}^2 - b_{21}^1 b_{12}^2)xy \\
&& + (b_{11}^{2}b_{22}^{2}-b_{12}^2 b_{21}^2)y^2
\end{eqnarray*} 
a polynomial that Mizner \cite{Mizner1} used in the study of codimension two CR structures and which
was mentioned earlier, independently, by Griffiths and Harris in the study of codimension two
subvarieties of complex projective space \cite{Griffiths-Harris1}.
Note that $P(x,y)$ is homogeneous of degree two in the variables $x$ and $y$. 
Let $A\in  GL(n, \mbox{{\bf C}})$ act on our bilinear form $B$.  Thus we have $B$ becoming
$(A^TB_1A,A^TB_2A)$.  Then the polynomial $P(x,y)$ transforms as follows:
$$\det (xA^TB_1A+ yA^TB_2A) = |\det(A)|^2 \det(xB_1 + yB_2)=   |\det(A)|^2 P(x,y).$$
By looking at this polynomial, we have effectively eliminated the influence of the 
$GL(n, \mbox{{\bf C}})$ action.  In other words, one method for generating invariants of the bilinear
form $B$ under the action of $GL(n, \mbox{{\bf C}})\times  GL(k, \mbox{{\bf C}})$ is to find invariants
of the polynomial $P(x,y)$ under the action of $ GL(k, \mbox{{\bf C}})$.  The action of $ GL(k, \mbox{{\bf C}})$
is just the standard change of basis on the variables $x$ and $y$.  For degree two homogeneous
polynomials in two variables, it is well know that the only invariant is the discriminant. (See
Chapter One of \cite{Olver}; recall, for the polynomial $Ax^2 + Bxy + Cy^2$, that the discriminant
is $B^2-4AC$.)  Thus for our bilinear form $B$, the invariant corresponding to the 
discriminant of the polynomial $\det(xB_1 + yB_2)$ will be
$$(b_{11}^1 b_{22}^2 + b_{22}^1 b_{11}^2
-b_{12}^1 b_{21}^2 - b_{21}^1 b_{12}^2)^2 - 4 (b_{11}^{1}b_{22}^{1}-b_{12}^1 b_{21}^1)
(b_{11}^{2}b_{22}^{2}-b_{12}^2 b_{21}^2),$$
which can be checked is $$(-1/2)v_{(23)(67)}\otimes w^{(23)}.$$

Again, we  have added this last part only to emphasis that there is geometry 
and meaning (though largely unexplored) behind the mechanical, almost crude, invariants that this 
paper generates.

\section{Conclusions}
There are many questions remaining.  One difficulty is in determining 
what a particular invariant or relation means.  This paper and 
\cite{GM}  just give
 lists, with no clue as to which have any type of important 
meaning, save for the invariants associated to the Mizner polynomial.
  (Of course, this is one of the difficulties in almost all of 
classical invariant theory).    Even 
harder is to determine when two vector-valued forms are equivalent.  
It is highly unlikely that the algebraic techniques in this paper will 
answer this question.


\begin{thebibliography}{99}

\bibitem{ACGH} Arbarello, E., Cornalba, M., Griffiths, P., and Harris, J., {\it Geometry
of Algebraic Curves}, Vol. I, Grundlehren der mathematischen Wissenschaften 267, Springer-Verlag,
New York, 1985.


\bibitem{Fulton} William Fulton and Joe Harris, {\it Representation Theory},
Springer-Verlag, New York, 1991.

\bibitem{DC} J.A. Dieudonne and J.B. Carrell, Invariant theory, old and
new, {\it Adv. Math.} 4, 1970, pp. 1-89.

\bibitem{Formanek}  E. Formanek, The invariants of $n\times n$ 
matrices, in {\it Invariant Theory} (S.S. Koh, Ed.), LNM 1268, 
Springer-Verlag, New York, 1987.

\bibitem{GM} Thomas Garrity and Robert Mizner, Invariants of
vector-valued bilinear and sesquilinear forms, {\it Linear Algebra
and its Applications}, Vol. 218 (1995), pp. 225-237.

\bibitem{Griffiths-Harris1} Griffiths, P. and Harris, J. Algebraic Geometry and Local
Differential Geometry, {\it Ann. Scient. Ec. Norm. Sup.} 12, 1979, pp 355-432.

\bibitem{Grossman} Zachary J. Grossman, Relations and Syzygies in 
Classical Invariant Theory for Vector-Valued Bilinear Forms, Williams 
College thesis, 1999.




\bibitem{Mizner1} Mizner, Robert I.,
CR structures of codimension $2$,
{\it J. Differential Geom.} 30 (1989), no. 1, 167--190. 


\bibitem{Popov}  V.L. Popov, {\it Groups, Generators, Syzygies,
and Orbits in Invariant Theory}, Translations of Mathematical
Monographs, vol. 100, AMS, Providence, 1992.

\bibitem{Procesi}  C. Procesi, The invariant theory of $n\times n$ 
matrices, {\it Advances in Mathematics}, vol. 19 (1976), pp. 306-381.

\bibitem{Olver} Peter J. Olver, {\it Classical Invariant Theory}, 
London Mathematical Society Student Texts 44, Cambridge University 
Press, 1999.


\bibitem{Sturmfels} Bernd Sturmfels, {\it Algorithms in Invariant Theory},
Springer-Verlag, New York, 1993.

\bibitem{Weyl} Hermann Weyl, {\it The Classical Groups},
Princeton University Press, Princeton, NJ, 1939.

\bibitem{Yokonuma}  Takeo Yokonuma, {\it Tensor Spaces and Exterior
Algebra}, Translations of Mathematical Monographs, vol. 108, AMS,
Providence, 1992.

\end{thebibliography}
\end{document}